\documentclass[11pt]{article}

\input{epsf}

\usepackage{amsfonts}

\newcommand{\qed}{\hbox{\rule{6pt}{6pt}}}

\newcommand{\Z}{\mathbb{Z}}

\newcommand{\LX}{ {\cal L} } 

\newtheorem{theorem}{Theorem}[section]
\newtheorem{corollary}[theorem]{Corollary}
\newtheorem{lemma}[theorem]{Lemma}

\newtheorem{example}[theorem]{Example}

\begin{document}

\title{Quandle Homology Theory and Cocycle Knot Invariants} 

\author{
J. Scott Carter\footnote{Supported in part by NSF Grant DMS \#9988107.}
\\University of South Alabama \\
Mobile, AL 36688 \\ carter@mathstat.usouthal.edu
\and 
Masahico Saito\footnote{Supported in part by NSF Grant DMS \#9988101.}
\\ University of South Florida
\\ Tampa, FL 33620  \\ saito@math.usf.edu
}

\maketitle

\vspace{10mm}

\begin{abstract}
This paper is a survey of  several 
papers in quandle homology theory and cocycle knot invariants 
that have been  published recently. 
Here we describe 
cocycle knot invariants 
that are defined in a state-sum form, 
quandle homology, and methods of constructing non-trivial cohomology classes.
\end{abstract}

\vspace{10mm}

\section{Prologue}

\begin{figure}
\begin{center}
\mbox{
\epsfxsize=4in
\epsfbox{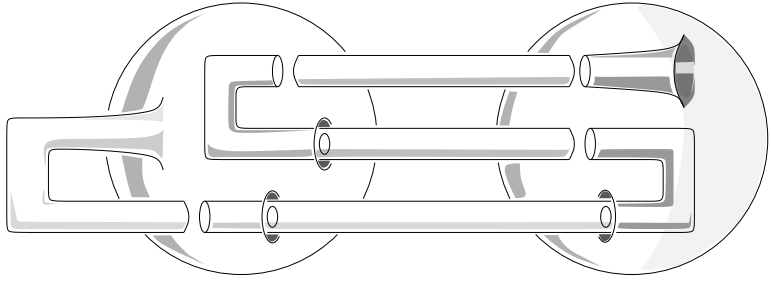} 
}
\end{center}
\caption{ Example 10 in ``Quick Trip'' }
\label{fmbroken} 
\end{figure}

We start with an example and its history. 
Figure~\ref{fmbroken} is an illustration of the knotted surface diagram
for an embedded $2$-sphere in the $4$-sphere, $S^4$. The $2$-sphere is
obtained by doubling a slice disk of the stevadore's knot. The diagram is
a {\it broken surface diagram} that is obtained from a generic projection
of the surface 
into $3$-space  by indicating over/under crossing
information in a way similar to the classical case. Specifically, the
portion of the surface that is closest to the hyperplane of projection is
depicted as an unbroken sheet while the sheet that is further away is
broken locally into two sheets. 
See \cite{CS:book} for details. 

\begin{figure}[h]
\begin{center}
\mbox{
\epsfxsize=4.5in
\epsfbox{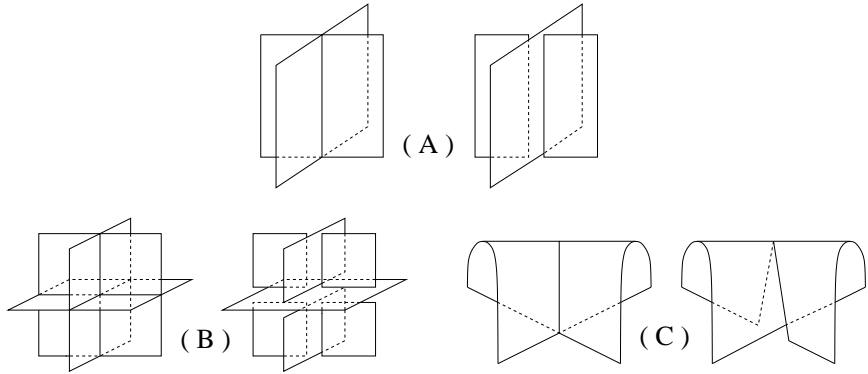} 
}
\end{center}
\caption{ Broken surface diagrams }
\label{broken} 
\end{figure}

Figure~\ref{broken} indicates the three local pictures at double, triple,
and branch points of the projection. 
A diagram can have branch and triple points in general, although the diagram
in Fig.~\ref{fmbroken} does not. 
At a triple point, 
we have a notion of
{\it top, middle and bottom} sheets. 
The adjectives describe the relative
proximity to the hyperplane into which the knotted surface has been
projected.

The sphere that is illustrated first appeared in
the manuscript~\cite{FM:57} by Fox and Milnor, 
  and later 
as Example 10 in Fox's ``Quick Trip'' \cite{FoxTrip},
described in a motion picture form.
This knotted sphere is not obtained by the spinning
construction~\cite{Artin}. 
This can be seen as follows.
The Alexander
polynomial of a spun knot  agrees with that of the underlying classical
knot since their fundamental groups are isomorphic.  
The first homology 
$H_1(\tilde{X})$ (called the knot module), 
of the infinite cyclic cover 
$\tilde{X}$  
of the complement $X$ of the 
sphere in $S^4$ 
in question, 
is   $\Z [ T, T^{-1}]/ (2-T)$   as a $\Lambda=\Z[T,T^{-1}]$-module,
thus the Alexander polynomial is not symmetric.

Fox's Example 11 can be recognized as the same sphere as Example 10 with
its orientation reversed. 
Its Alexander polynomial is $(1-2T)$.
Thus the
sphere illustrated in Fig.~\ref{fmbroken} is 
{\it non-invertible}: 
It is
not ambiently isotopic to the same surface with its orientation reversed.

Example 12 of ``Quick Trip'' 
has as its knot module $\Lambda/(2-T,1-2T)$.
It is obtained from 
the previous  two examples 
by combining  some of their portions. 
The
fact that this ideal is not principal also illustrates the difference
between classical knot theory and knotted surfaces. 
Note that the argument of asymmetric ideals no longer 
applies to Example 12.  
It is also interesting to
note that this example is in fact the $2$-twist spun trefoil~\cite{Lith}, 
although
Zeeman's  twist spin construction appeared later
in 1965 \cite{Ze65}.

Hillman \cite{Hillman} showed that this knotted sphere was 
non-invertible 
using the Farber-Levine pairing.
Ruberman \cite{Ruber}  
used Casson-Gordon invariants  
to prove the same result,
with other new examples of non-invertible knotted spheres. 
Neither technique applies
directly to the same knot with a trivial $1$-handle attached.
Kawauchi~\cite{Kawa86,Kawa90a}  
has generalized the Farber-Levine pairing to higher genus surfaces, 
showing that such a 
torus is also non-invertible. 
The method we survey in this article shows this fact \cite{CJKLS} using
an invariant defined in a state-sum form from 
quandle cohomology theory,
called the cocycle knot invariant. 
The cocycle knot invariant has also been used to prove
new geometric results \cite{SatShima}.

We asked Ruberman if he had proved 
non-invertibility of the $2$-twist-spun trefoil 
 on his first excursion to the
Georgia Topology Conference in 1982. (Incidentally, the first named author
also had his topology debute at GTC1982. The second named author debued
 at GTC1990.)  
Ruberman told us that the era was correct, although 
he did not present the result then. 
His dissertation, however, was inspired by the paper by Sumners~\cite{Sumners},
which showed, in  particular, that  any $2$-sphere in $4$-space
that contains the Stevedore's knot as a cross-section is knotted,
such as the above examples in ``Quick Trip.''

In Section~\ref{quanss} below, we will give the definition of 
the cocycle invariant 
for classical knots and for knotted surfaces in $4$-space. 
Our motivation came from 
the Jones polynomial and quantum 
invariants of $3$-manifolds. 
A common feature of the quantum invariants is the state-sum definition,
and it has been asked since their discovery whether such invariants 
exist in higher dimensions (see \cite{CKS,CKY} 
 for such  attempts). 
We briefly review the  state-sum definition of 
Jones polynomial and a related invariant 
for  triangulated 
$3$-manifolds --- the Dijkgraaf-Witten invariant.

\begin{figure}
\begin{center}
\mbox{
\epsfxsize=3in
\epsfbox{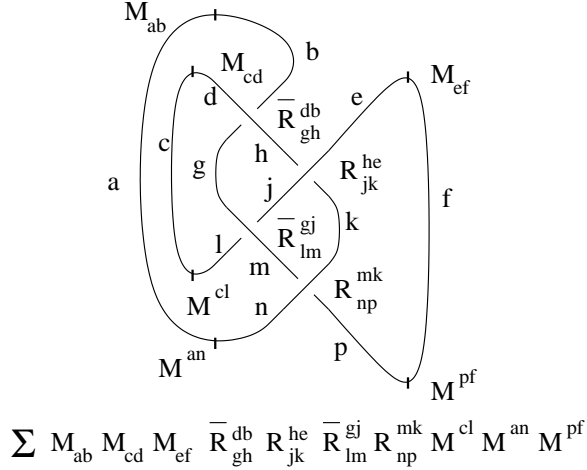} 
}
\end{center}
\caption{ Bracket polynomial of knots }
\label{bracket} 
\end{figure}

\subsection*{The Bracket Model} 

The bracket polynomial of a classical knot or link is obtained as follows.
The knot is projected generically into the plane, and a height function on
the plane 
is chosen.  
Let an index set $S=\{1,2\}$ 
(in general a finite set),
whose elements are called spins, 
be given and fixed. 
Let ${\cal A}$ denote the set of arcs obtained from the given 
knot diagram by deleting local maxima, minima, and crossing points. 
The {\it coloring} ${\cal C}$  is a map ${\cal C}: {\cal A} \rightarrow S$.

{\it Boltzmann weights} $B(\tau, {\cal C})$ 
are assigned at minima, maxima, and crossings
as follows: 
Local minima are assigned 
$M^{ab}$, 
local maxima are assigned $M_{ab}$, 
 crossings are assigned $R^{ab}_{cd}$ 
if the over crossing arc has positive slope, or 
$\overline{R}^{ab}_{cd}$ 
if the over-crossing arc has negative slope, 
where each weight is defined with a variable $A$ and $i=\sqrt{-1}$ by
$$ M^{ab}=M_{ab}=
\left\{ \begin{array}{ll} iA & \mbox{if \quad $a=0$, $b=1$} \\
(iA)^{-1} &  \mbox{if \quad $a=1$, $b=0$} \\
0  &  \mbox{otherwise}
\end{array} \right. , 
$$

\vspace{-4mm}

\begin{eqnarray*}
 R^{ab}_{cd} & = & A \delta^a_c \delta^b_d + A^{-1} M^{ab} M_{cd}  ,
\\
\overline{R}^{ab}_{cd} &=&  A^{-1} \delta^a_c \delta^b_d + A M^{ab} M_{cd} .
\end{eqnarray*} 
Here, $\delta$ denotes
Kronecker's delta. 
 The bracket polynomial, as a
polynomial in $A$, is defined by the {\it state-sum} 
$$ \langle K \rangle = \sum_{\cal C} \prod_{\tau}  B(\tau, {\cal C} ),$$
where the product is taken over all crossings, and the sum is taken 
over all colorings.
The
Jones polynomial is obtained from the bracket by normalizing and
substituting. 
Spefically,  the quantity 
${\mathcal L}_K (A) = (-A)^{-3w} \langle K \rangle$
is a knot invariant, where the exponent $w$ is the {\it
writhe} of the diagram $K$, and  
$V(t)= {\mathcal L}_K (t^{-1/4})$ is the Jones polynomial.
In Fig.~\ref{bracket},
a colored knot diagram and its Boltzmann weights are depicted.  
For a given coloring (denoted by lower case 
letters $a$ through $p$ excluding $i$ ) 
on arcs, the product of the Boltzmann weights are given at the bottom
of the figure. The sum is taken over all colorings,
see~\cite{K&P}
for details. 

\subsection*{The Dijkgraaf-Witten Invariant}

\begin{figure}
\begin{center}
\mbox{
\epsfxsize=1.3in
\epsfbox{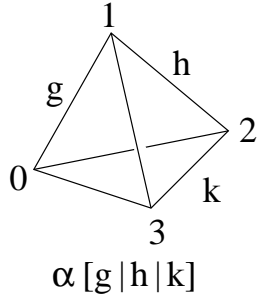}
}
\end{center}
\caption{A $3$-cocycle assigned at a triangle}
\label{dw}
\end{figure}

Similar state-sum invariants were defined 
for $3$-manifolds 
in \cite{DW} using group cocycles and the state-sum 
concept as follows. A combinatorial
definition for Chern-Simons invariants 
with finite gauge groups  was given 
using $3$-cocycles of 
group cohomology.
We follow Wakui's description,
see \cite{Wakui} for more detailed treatments.
Let $T$ be a triangulation of an oriented closed
$3$-manifold $M$,
with $a$ vertices and $n$ tetrahedra.
Give an ordering to the set of vertices. Let $G$ be a finite group.
Let ${\cal C} : $ \mbox{ $\{$ oriented edges $\} \rightarrow G$ } 
be a map such 
that
(1) for any triangle with vertices $v_0, v_1, v_2$ of $T$,
${\cal C}(\langle v_0, v_2 \rangle )={\cal C}( \langle v_1, v_2 \rangle )
{\cal C} ( \langle v_0, v_1 \rangle )$,
where $ \langle v_i, v_j \rangle $ denotes the oriented edge
with endpoints $v_1$ and $v_2$, 
and
(2) ${\cal C}(-e) = {\cal C}(e)^{-1}$.
Such a map ${\cal C}$ is called a (group) coloring. 
Let $\alpha :G \times G \times G \rightarrow A$,
$(g,h,k) \mapsto \alpha [g|h|k] \in A$, be a $3$-cocycle
with   values  in a 
multiplicative abelian group $A$, $\alpha \in Z^3(G;A)$. 
The $3$-cocycle condition is
written as 
$$\alpha [h|k|l] \alpha [gh|k|l]^{-1} \alpha [g|hk|l]
\alpha [g|h|kl]^{-1} \alpha [g|h|k] =1. $$
Then the Dijkgraaf-Witten invariant is defined by
$$ Z_M = \frac{1}{|G|^a} \sum_{{\cal C}} \prod_{i=1}^{n} W( \sigma, {\cal C} )^{\epsilon_i}. $$ 
Here $a$ denotes the number of the vertices of the
given triangulation,                     
$W(\sigma, {\cal C})= \alpha [g|h|k] $ where
${\cal C}( \langle v_0, v_1 \rangle )=g$, ${\cal C}( \langle v_1, v_2 \rangle  )=h$,
 ${\cal C}( \langle v_2, v_3 \rangle )=k$,
for the tetrahedron $\sigma = |v_0 v_1 v_2 v_3|$ with
the ordering $v_0<v_1<v_2<v_3$,
and  
$\epsilon=\pm 1 $ according to 
whether
or not the orientation
of $\sigma$ with respect to the vertex ordering matches the
orientation of $M$,
see Fig.~\ref{dw}. 

\section{Quandles and Quandle Colorings}

In this section we 
define quandles, quandle 
colorings, and illustrate that counting  quandle colorings can be
formulated as a state-sum. This definition will help motivate the definition
of the quandle cocycle invariants that we will define in Section~\ref{quanss}.

A {\it quandle}, $X$, is a set with a binary operation 
$(a, b) \mapsto a * b$
such that

(I) For any $a \in X$,
$a* a =a$.

(II) For any $a,b \in X$, there is a unique $c \in X$ such that 
$a= c*b$.

(III) 
For any $a,b,c \in X$, we have
$ (a*b)*c=(a*c)*(b*c). $

A {\it rack} is a set with a binary operation that satisfies 
(II) and (III).

\begin{figure}[h]
\begin{center}
\mbox{
\epsfxsize=4in
\epsfbox{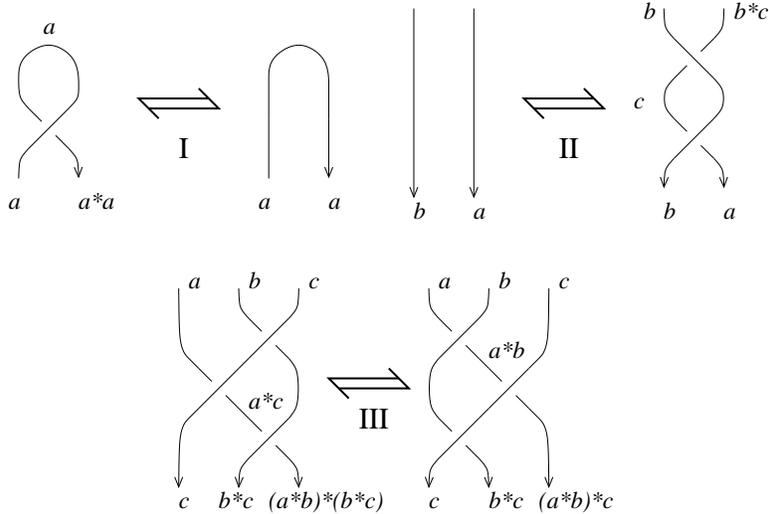}
}
\end{center}
\caption{Reidemeister moves and quandle conditions}
\label{Rmoves}
\end{figure}

Racks and quandles have been studied in, for example, 
\cite{Br88,FR,Joyce82,K&P,Matveev}.
The axioms for a quandle correspond respectively to the 
Reidemeister moves of type I, II, and III 
(see Fig.~\ref{Rmoves} and 
\cite{FR,K&P}, for example). 
A function $f: X \rightarrow  Y$ between quandles
or racks  is a {\it homomorphism}
if $f(a \ast b) = f(a) * f(b)$ 
for any $a, b \in X$. 
The following are typical examples of quandles.

\begin{itemize}
\item
A group $X=G$ with
$n$-fold conjugation
as the quandle operation: $a*b=b^{-n} a b^n$.
\item
Any set $X$ with the operation $x*y=x$ for any $x,y \in X$ is
a quandle called the {\it trivial} quandle.
The trivial quandle of $n$ elements is denoted by $T_n$.
\item
Let $n$ be a positive integer.
For elements  $i, j \in \{ 0, 1, \ldots , n-1 \}$, define
$i\ast j \equiv 2j-i \pmod{n}$.
Then $\ast$ defines a quandle
structure  called the {\it dihedral quandle},
  $R_n$.
This set can be identified with  the
set of reflections of a regular $n$-gon
  with conjugation
as the quandle operation.
\item
Any $\Lambda (={\Z }[T, T^{-1}])$-module $M$
is a quandle with
$a*b=Ta+(1-T)b$, $a,b \in M$, called an {\it  Alexander  quandle}.
Furthermore for a positive integer
$n$, a {\it mod-$n$ Alexander  quandle}
${\Z }_n[T, T^{-1}]/(h(T))$
is a quandle
for
a Laurent polynomial $h(T)$.
It 
is finite if the coefficients of the
highest and lowest degree terms
of $h$   are units in $\Z_n$. 
\end{itemize}

\begin{figure}
\begin{center}
\mbox{
\epsfxsize=2.5in
\epsfbox{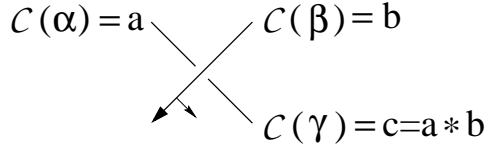} 
}
\end{center}
\caption{ Quandle relation at a crossing  }
\label{qcolor} 
\end{figure}

Let $X$ be a fixed quandle.
Let $K$ be a given oriented classical knot or link diagram,
and let ${\cal R}$ be the set of (over-)arcs. 
The normals are given in such a way that (tangent, 
normal) matches
the orientation of the plane, see Fig.~\ref{qcolor}. 
A (quandle) {\it coloring} ${\cal C}$ is a map 
${\cal C} : {\cal R} \rightarrow X$ such that at every crossing,
the relation depicted in Fig.~\ref{qcolor} holds. 
More specifically, let $\beta$ be the over-arc at a crossing,
and $\alpha$, $\beta$ be under-arcs such that the normal of the over-arc
points from $\alpha$ to $\beta$.
Then it is required that ${\cal C}(\gamma)={\cal C}(\alpha)*{\cal C}(\beta)$.

Alternately, a coloring can be described as a 
quandle homomorphism as follows.
Classical knots have fundamental quandles that are defined via generators
and relations. The theory of quandle presentations is given a complete
treatment in \cite{FR}. Specifically, the generators of the fundamental
quandle correspond to the arcs in a diagram. 
The quandle relation
$a*b=c$  holds where $a$ is the generator that corresponds to the underarc
away from which the normal to the over arc points, $b$ is the generator
that corresponds to the overarc, and $c$ corresponds to the 
underarc towards which the transversal's normal points,
see Fig.~\ref{qcolor}. 
A  coloring  of a classical knot diagram by 
a quandle $X$ gives rise to 
a quandle homomorphism from the fundamental quandle to the quandle $X$.

The number $\mbox{Col}_X(K)$  of colorings of a knot diagram 
$K$ by a fixed finite quandle $X$ is a knot invariant, and 
has  a  description as a state-sum 
as follows. 

For a finite quandle $X$, consider 
the set of maps $\{ {\cal D} : {\cal R} \rightarrow  X \}$ 
(without the requirement of a quandle coloring). 
For a given such a map ${\cal D}$, 
define the Boltzmann weight at a crossing $\tau$, with over-arc
$\beta$ whose normal points from the under-arc $\alpha$ to the under-arc
$\gamma$,  by
$$B(\tau, {\cal D} )= 
 \left\{ \begin{array}{ll} 1 & {\mbox{\rm if}} \ \ 
 {\cal D}(\alpha)*  {\cal D}(\beta)=  {\cal D}(\gamma)  \\ 
0 & {\mbox{\rm otherwise}}. \end{array} \right. $$
Then the number of quandle  colorings is written by a state-sum
$ \mbox{Col}_X(K)=\sum_{{\cal D} } \prod_{\tau} B(\tau, {\cal D})$.
We could also use  colorings similar to those used in the bracket, 
or we could write 
$ \mbox{Col}_X(K)=\sum_{{\cal C} } \prod_{\tau} B_1(\tau, {\cal C})$, 
where ${\cal C} $ ranges over only quandle 
colorings ${\cal C}$,
and $B_1(\tau, {\cal C})\equiv 1$ is a constant function. 
Either way, it is natural to ask whether we can modify 
the weights $1$ to a general function.

Fox's $n$-coloring is a quandle coloring by the dihedral quandle $R_n$.
The classical result that a knot is 
non-trivially 
Fox $p$-colorable
if and only if 
$p| \Delta(-1)$ 
(where $\Delta(t)$ 
denotes the Alexander polynomial)
has been generalized by Inoue~\cite{Inoue} to the following:

Let $\Delta_K^{(i)}(T)$ denote the greatest common divisor of 
all $(n-i-1)$ minor determinants of the presentation matrix 
for the knot module obtained via the Fox calculus.

\begin{theorem}{\bf \cite{Inoue}} Let $p$ be a prime number, 
$J$ an ideal of the ring $\Lambda_p=\Z_p[T,T^{-1}]$ 
and let $Q(K)$ denote a knot quandle. For each 
$i\ge 0$, put  
$e_i(T)=\Delta_{K}^{(i)}(T)/\Delta_{K}^{(i+1)}(T).$
Then the number of all quandle homomorphisms of the knot quandle $Q(K)$
 to the Alexander quandle 
$\Lambda_p/J$ is equal to the cardinality of the module 
$\Lambda_p/J \oplus \oplus_{i=0}^{n-2} \{ \Lambda_p/(e_i(T),J)\}$
\end{theorem}

\begin{figure}
\begin{center}
\mbox{
\epsfxsize=3in 
\epsfbox{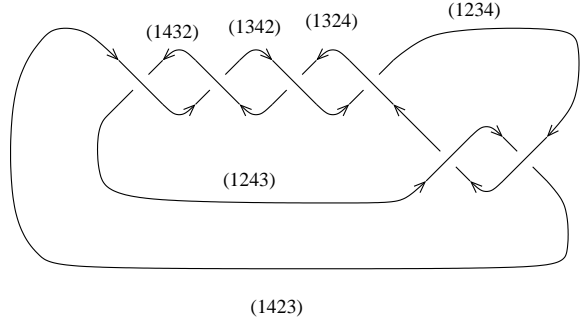} 
}
\end{center}
\caption{ A coloring of $6_1$ by   $QS_6$ }
\label{six1} 
\end{figure}

\begin{example}{\rm The Alexander quandle $S_4=\Z_2[T,T^{-1}]/(T^2+T+1)$
has four elements that are 
represented as $0,1,T,$ and $T+1$.
This quandle colors both 
 the trefoil ($3_1$) and  the figure $8$ knot ($4_1$) 
as one can easily see directly or by considering
the mod-$2$ reduction of the Alexander polynomials. 
In either case, the order of 
$S_4$ is $4$ but the determinants
are $3$ and
$5$, for $3_1$ and $4_1$, respectively. 
Thus quandle colorings are more general than Fox colorings.
}\end{example}

\begin{example}{\rm The quandle $QS_6$ consists of the 
$4$-cycles
$a=(1234)$, $A=(1432)$, $b=(1243)$, 
$B=(1342)$, $c=(1324)$, and $C=(1423)$
with group conjugation as the quandle operation. 
Figure~\ref{six1} illustrates a coloring of the knot $6_1$ by $QS_6$. 
This quandle has $R_3$ as a quotient quandle. 
The map $f(a)=f(A)=0$, $f(b)=f(B)=1$, and $f(c)=f(C)=2$ 
is a quandle homomorphism. 
The equalizers ($E_y=\{x: f(x)=f(y)\}$)
are all the two element trivial quandle.
Recently,
Angela Harris has shown that $QS_6$ is not an Alexander quandle
 of the form $\Lambda_n/(h)$ where $h$ is a polynomial.

}\end{example}

\section{Quandle Homology and Cohomology Theories} 

In this section, we present  twisted quandle homology,
which was discussed in \cite{CENS}, and specialize it to the untwisted
theory subsequently. 
Originally, rack homology and homotopy theory were defined and 
studied in \cite{FRS}, and a modification to quandle homology theory 
was given in \cite{CJKLS} to define a knot invariant in a state-sum form.
Then they were  generalized to 
a  twisted theory in \cite{CENS}.
Computations are found in \cite{CJKS1,betti}  and also in
\cite{LN,Mochi} by other authors.

Let $\Lambda=\Z [T, T^{-1}]$, and let 
$C_n^{\rm TR}(X)= C_n^{\rm TR}(X; \Lambda)$ be the free 
module over $\Lambda$   generated by
$n$-tuples $(x_1, \dots, x_n)$ of elements of a quandle $X$. 
Define a homomorphism 
$\partial = \partial^T_{n}: C_{n}^{\rm TR}(X) \to  C_{n-1}^{\rm TR}(X) $ 
by \begin{eqnarray*} 
\lefteqn{
\partial^T_{n}(x_1, x_2, \dots, x_n) } \nonumber \\ && =
\sum_{i=1}^{n} (-1)^{i}\left[ T (x_1, x_2, \dots,
x_{i-1}, x_{i+1},\dots, x_n) \right.
\nonumber \\
&&
- \left. (x_1 \ast x_i, x_2 \ast x_i, \dots, x_{i-1}\ast x_i, x_{i+1}, 
\dots, x_n) \right]
\end{eqnarray*}
for $n \geq 2$ 
and $\partial^T_n=0$ for 
$n \leq 1$. 
We regard that the $i=1$ terms contribute $(1-T) (x_2, \ldots, x_n)$. 
Then
$C_\ast^{\rm TR}(X) = \{C_n^{\rm TR}(X), \partial^T_n \}$
is a chain complex.
For any $\Lambda$-module $A$, let 
$C_\ast^{\rm TR}(X;A)= 
\{C_n^{\rm TR}(X) \otimes _{\Lambda} A , \partial^T_n \}$
be the induced chain complex, where the induced boundary operator is 
represented by the same notation.
Let $C^n_{\rm TR}(X;A)=\mbox{Hom}_{\Lambda} ( C_n^{\rm TQ}(X), A)$
and define the coboundary operator 
$\delta=\delta^n_{\rm TR}  : C^n_{\rm TR}(X;A) \to C^{n+1}_{\rm TR}(X;A)$
by 
$(\delta f )(c)= (-1)^n f( \partial c) $ 
for any $c \in  C_n^{\rm TQ}(X)$
and $f \in C^n_{\rm TR}(X;A)$. Then 
$C^\ast_{\rm TR}(X;A)= 
\{C^n_{\rm TR}(X;A) , \delta^n_{\rm TR}  \}$
is a cochain complex.
The $n\/$-th homology and cohomology groups of 
these complexes 
are called 
{\it twisted rack homology group\/} and {\it cohomology group\/},
and are denoted by $H_n^{\rm TR}(X; A)$ and $H^n_{\rm TR}(X; A)$,
respectively.

Let $C_n^{\rm TD}(X;A)$ be the subset of $C_n^{\rm TR}(X;A)$ generated
by $n$-tuples $(x_1, \dots, x_n)$
with $x_{i}=x_{i+1}$ for some $i \in \{1, \dots,n-1\}$ if $n \geq 2$;
otherwise let $C_n^{\rm TD}(X;A)=0$. If $X$ is a quandle, then
$\partial^T_n(C_n^{\rm TD}(X;A)) \subset C_{n-1}^{\rm TD}(X;A)$ and
$C_\ast^{\rm TD}(X;A) = \{ C_n^{\rm TD}(X;A), \partial^T_n \}$ is 
a sub-complex of
$C_\ast^{\rm TR}(X;A)$. 
Similar subcomplexes 
$C^\ast_{\rm TD}(X;A) = \{ C^n_{\rm TD}(X;A), \delta_T^n \}$
are defined for cochain complexes.
The $n\/$-th homology and cohomology groups of 
these complexes  
are called 
{\it twisted degeneracy  homology group\/} and {\it cohomology group\/},
and are denoted by $H_n^{\rm TD}(X; A)$ and $H^n_{\rm TD}(X; A)$,
respectively.

Put $C_n^{\rm TQ}(X;A) = C_n^{\rm TR}(X;A)/ C_n^{\rm TD}(X;A)$ and 
$C_\ast^{\rm TQ}(X;A) = \{ C_n^{\rm TQ}(X;A), {\partial  }^T_n \}$, 
where all the induced boundary 
and coboundary 
operators are denoted by 
$\partial =  {\partial  }^T_n $
and $\delta={\delta}_T^n $, respectively.   
A cochain complex 
$C^\ast_{\rm TQ}(X;A) = \{ C^n_{\rm TQ}(X;A), {\delta }_T^n \}$
is similarly defined.
The $n\/$-th homology and cohomology groups of 
these complexes 
are called 
{\it twisted   homology group\/} and {\it cohomology group\/},
and are denoted by
\begin{eqnarray*}
H_n^{\rm TQ}(X; A) 
= H_{n}(C_\ast^{\rm TQ}(X;A)), \quad
H^n_{\rm TQ}(X; A) 
= H^{n}(C^\ast_{\rm TQ}(X;A)). \end{eqnarray*}
The groups of (co)cycles and (co)boundaries are denoted 
using similar notations.

\begin{example} \label{2cocyex} {\rm
The $1$-cocycle condition is written for $\eta \in Z^1_{ TQ}(X;A)$  as 
$$ -T \eta(x_2) + T \eta(x_1) + \eta(x_2) - \eta(x_1*x_2)=0, \quad  \mbox{or}$$

\vspace{-6mm} 

$$ T \eta(x_1) + (1-T)\eta(x_2)  = \eta(x_1*x_2) . $$
Note that this means that $\eta: X \rightarrow A$ is a quandle homomorphism.

The $2$-cocycle condition is written for $\phi \in Z^2_{ TQ}(X;A)$  as 
\begin{eqnarray*}
\lefteqn{ T [ - \phi (x_2, x_3) + \phi (x_1, x_3) - \phi (x_1, x_2) ] } \\
& + & [ \phi (x_2, x_3) -  \phi (x_1 * x_2, x_3) 
+ \phi (x_1 * x_3, x_2 * x_3) ] =0 \quad \mbox{or}
\end{eqnarray*}

\vspace{-8mm}

\begin{eqnarray*}
\lefteqn{ T  \phi (x_1, x_2) +   \phi (x_1 * x_2, x_3) } \\
& = & T  \phi (x_1, x_3) + (1-T) \phi (x_2, x_3)
+  \phi (x_1 * x_3, x_2 * x_3).
\end{eqnarray*}

The  
geometric meaning of this condition will become clear 
in Section~\ref{quanss}. 

} \end{example}
The original untwisted quandle homology is described as 
a specification of $T=1$. 
Specifically, in the definition of the boundary homomorphism 
$\partial ^T_n$, set $T=1$, and define all the cycle, boundary, homology groups
similarly. Then use Hom$( \ - \ ; A)$ 
to define cohomology theory. 
Thus we assume that the coefficients  $A$ simply form   an abelian group.
We obtain 
degenerate, rack, and quandle homology groups denoted by 
$H^{\rm W}_*(X;A)$ for $W=D,R,Q$, respectively. 
Similarly, $H_{\rm W}^n(X;A)$ denotes the corresponding cohomology groups. 
The cohomology theory $H_{\rm R}^*$ was defined in \cite{FRS}. 
It was seen in
\cite{betti} 
that the short exact sequence:
\begin{eqnarray*}
0 \to C_n^{\rm D}(X) \stackrel{i}{\to} C_n^{\rm R}(X) \stackrel{j}{\to}
C_n^{\rm Q}(X) \to 0 
\end{eqnarray*}
gives rise to
the following homology long exact sequence:
\begin{eqnarray*}
\cdots \stackrel{\partial_\ast}{\to} H_n^{\rm D}(X;A) \stackrel{i_\ast}{\to}
H_n^{\rm R}(X;A)
\stackrel{j_\ast}{\to} H_n^{\rm Q}(X;A)
\stackrel{\partial_\ast}{\to} H_{n-1}^{\rm D}(X;A) \to \cdots
\end{eqnarray*}
and it was shown in  \cite{CKS:geo} by geometric arguments  
that the sequence splits in low dimensions.
This result was improved upon in \cite{LN} 
by Litherland and Nelson 
where they showed the following:

\begin{theorem}{\bf \cite{LN}} The above
long exact sequence splits
into short exact sequences 
$$ 0 {\to}  H_n^{\rm D}(X;A) {\to}
H_n^{\rm R}(X;A)
{\to} H_n^{\rm Q}(X;A)
{\to} 0 . $$ 
\end{theorem}

In fact, they 
construct a projection 
$p:C_n^{\rm R}(X) {\to} C_n^{\rm D}(X)$ 
thereby splitting the short exact sequence of chain complexes.

\section{Cocycle Knot Invariants} 
\label{quanss}

\subsection*{Untwisted Cocycle Invariants}

Let $K$ be a classical knot or link diagram. Let a finite 
quandle $X$, and an (untwisted) 
quandle $2$-cocycle $\phi \in Z^2_{\rm Q}(X;A)$ be given.
A {\it (Boltzmann) weight}, $B(\tau, {\cal C})$ 
(that depends on $\phi$),
at a  crossing $\tau$ is defined as follows.
Let  ${\cal  C}$ 
denote a coloring ${\cal  C}: {\cal R} \rightarrow X$. 
Let $\beta$ be the over-arc at $\tau$, and $\alpha$, $\gamma$ be 
under-arcs such that the normal to $\beta$ points from $\alpha$ to $\gamma$,
see Fig.~\ref{qcolor}. 
Let $x={\cal C}(\alpha)$ and $y={\cal C}(\beta)$. 
Then define $B(\tau, {\cal C})= \phi(x,y)^{\epsilon (\tau)}$,
where $\epsilon (\tau)= 1$ or $-1$, if  
(the sign of) the crossing  $\tau$ 
is positive or negative, respectively.
By convention, the crossing in Fig.~\ref{qcolor} is positive if 
the orientation of the under-arc points downward.

\begin{figure}
\begin{center}
\mbox{
\epsfxsize=4in
\epsfbox{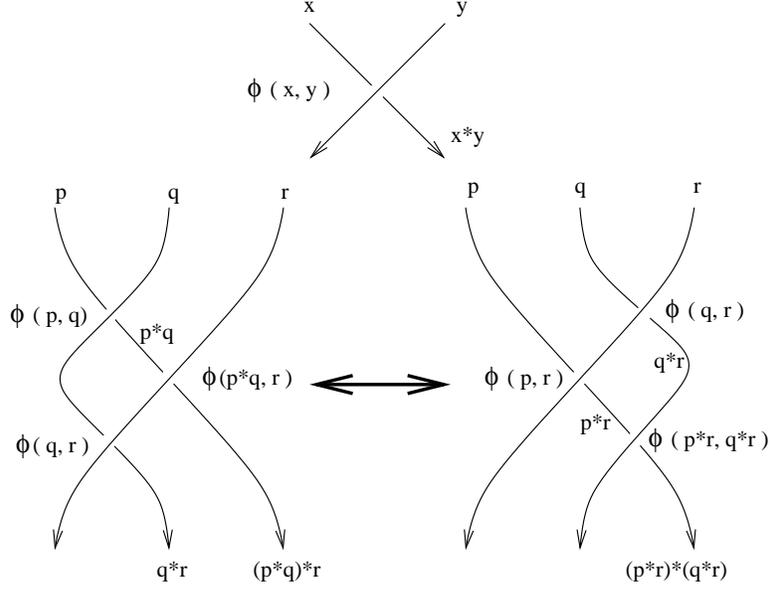} 
}
\end{center}
\caption{ The untwisted $2$-cocycle condition and type III move  }
\label{2cocy} 
\end{figure}

The {\it (quandle) cocycle knot invariant} is defined by 
the state-sum expression  
$$
\Phi (K) = \sum_{{\cal C}}  \prod_{\tau}  B( \tau, {\cal C}). 
$$
The product is taken over all crossings of the given diagram $K$,
and the sum is taken over all possible colorings.
The values of the partition function 
are  taken to be in  the group ring ${\bf Z}[A]$ where $A$ is the coefficient 
group  written multiplicatively. 
The state-sum depends on the choice of $2$-cocycle $\phi$. 
This is proved \cite{CJKLS} to be a knot invariant.  
Figure~\ref{2cocy} shows the invariance of the state-sum under
the Reidemeister type III move.
The sums of cocycles, equated before and  after the move, 
is the $2$-cocycle condition given in Example~\ref{2cocyex} 
with the evaluation $T=1$. 

The following variations have been considered.
\begin{itemize}
\item
Lopes~\cite{Lopes}  
observed that the family
$\{ \prod_{\tau} B(\tau, {\cal C} ) \}_{\cal C} $ 
is a knot invariant, without taking summation.
In particular, infinite quandles can be used for coloring
in this case.

\item
For a link $L=K_1 \cup \ldots \cup K_n$, let ${\cal T}_i$, 
$i=1, \ldots, n$, be the set of crossings at which the under-arcs
belong to the component $K_i$. Then it was observed \cite{CENS} 
that $\{ \sum_{\cal C} \prod_{ \tau \in {\cal T}_i }  B(\tau, {\cal C} ) \}_{i=1}^n $ is a link invariant, strictly stronger than the single state-sum.

\end{itemize}

\subsection*{Twisted Cocycle Invariants}

Let $K$ be an oriented knot diagram with normals. 
The (underlying) diagram divides the plane into regions.
Take an arc $\ell$ from  the region at infinity to 
a region 
$H$ such that $\ell$  intersects the arcs (missing crossings) of 
the diagram transversely in  finitely many points.
The {\it Alexander numbering} $\LX(H)$ of a region $H$ is the number 
of such intersections counted with signs. 
This does not depend on the choice of an arc $\ell$.

\begin{figure}
\begin{center}
\mbox{
\epsfxsize=2.5in
\epsfbox{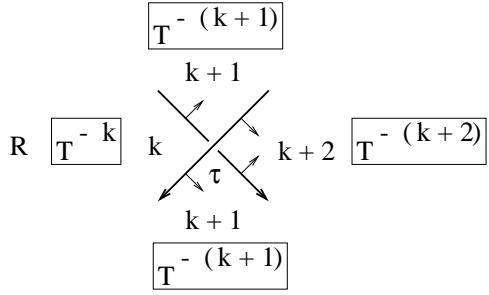} 
}
\end{center}
\caption{Alexander numbering of a crossing }
\label{LXcrossing} 
\end{figure}

Let $\tau$ be a crossing.
There are four regions near $\tau$, 
and  the unique region
from which normals of over- and under-arcs point  
is called the {\it source region} of $\tau$.
The {\it Alexander numbering} $\LX(\tau)$ of a crossing $\tau$ 
is defined to be 
$\LX(R) $ 
where $R$ is the source region of $\tau$. 
Compare with \cite{CKamS:LX}.
In other words, $\LX(\tau)$ is the number of intersections, 
counted with signs,
between an arc $\ell$ from the region at infinity to
$\tau$ approaching from  the source region
of $\tau$.
In  Fig.~\ref{LXcrossing}, the source region $R$ is the left-most region,
and the Alexander numbering of $R$ is $k$,
and so is the Alexander numbering of the crossing $\tau$.

Let a classical knot (or link) diagram $K$,  a finite quandle $X$,
a finite Alexander quandle $A$  
be given.  
A coloring 
of $K$ by $X$  also is given and 
is denoted by ${\cal C}$. 
A {\it twisted (Boltzmann) weight}, $B_T(\tau, {\cal C})$,
at a  crossing $\tau$ is defined as follows.
Let  ${\cal  C}$ 
denote a coloring. Let $\beta$ be the over-arc at $\tau$, and $\alpha$, $\gamma$ be 
under-arcs such that the 
normal
to $\beta$ points from  $\alpha$ to $\gamma$.
Let $x={\cal C}(\alpha)$ and $y={\cal C}(\beta)$. 
Pick a 
twisted 
 quandle  $2$-cocycle  $\phi \in  Z^2_{\rm TQ}(X; A)$.
Then define 
$B_T(\tau, {\cal C})= [\phi(x,y)^{\epsilon (\tau)} ]^{T^{-\LX(\tau)}}$,
where  $\epsilon (\tau)= 1$ or $-1$, if  the sign of $\tau$ 
is positive or negative, respectively.
Here, we use the multiplicative notation of elements of $A$, so that 
$\phi(x,y)^{-1}$ denotes the inverse of $\phi(x,y)$. 
Recall that $A$ admits an action by $\Z=\{ T^n \}$, 
and for $a \in A$, the action of $T$ on $a$ is denoted by $a^T$. 
To specify  the action by ${T^{-\LX(\tau)}}$ in the figures,
each region $R$ with 
Alexander numbering
$\LX(R)=k$ is labeled  by the power $T^{-k}$ framed with a square,
as depicted in Fig.~\ref{LXcrossing}.

The {\it state-sum}, or a {\it partition function},
is the expression 
$$
\Phi_T (K) = \sum_{{\cal C}}  \prod_{\tau}  B_T( \tau, {\cal C}).
$$  
The product is taken over all crossings of the given diagram,
and the sum is taken over all possible colorings.
The value of the weight  $B_T( \tau, {\cal C})$
is in the coefficient group $A$ written multiplicatively. 
Hence the value of the state-sum is in  the group ring ${\bf Z}[A]$.

\begin{figure}[h]
\begin{center}
\mbox{
\epsfxsize=4in
\epsfbox{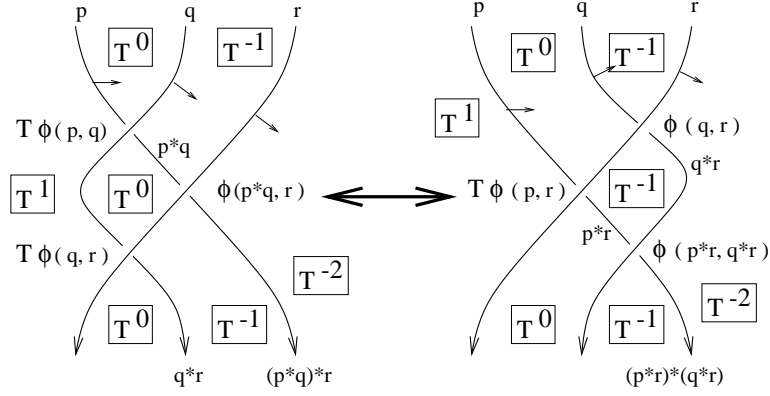} 
}
\end{center}
\caption{ The twisted $2$-cocycle condition and type III move }
\label{2tcocy} 
\end{figure}

It was proved in  \cite{CES} that $\Phi_T (K)$ is a knot invariant,
called the {\it (quandle) twisted cocycle invariant}. 
Figure~\ref{2tcocy} depicts the invariance under the type III move,
where the left-most region is assumed to have Alexander numbering $-1$. 
The sum over all these cocycles, equated before and after the move,
gives the $2$-cocycle condition written in Example~\ref{2cocyex}.

\subsection*{Cocycle 
Invariants for Knotted Surfaces}

The state-sum invariant is defined in an analogous way for 
oriented knotted surfaces
in $4$-space using their projections and diagrams in $3$-space. 
Specifically, the above steps can be repeated as follows,
for a fixed finite quandle $X$ and a knotted surface diagram $K$.

\begin{figure}[h]
\begin{center}
\mbox{
\epsfxsize=3.5in
\epsfbox{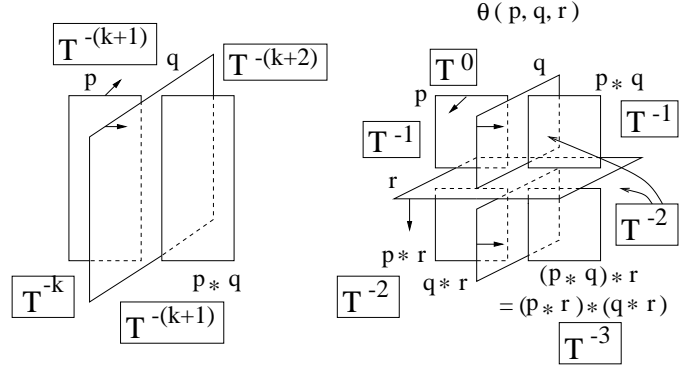} 
}
\end{center}
\caption{Colors at double curves and $3$-cocycle at a triple point }
\label{triplepoint} 
\end{figure}

\begin{itemize}
\item
The diagrams consist of double curves and isolated branch and triple points
\cite{CS:book}. Along the double curves, 
the 
coloring rule is defined using normals
in the same way as classical case, as depicted in
the left of Fig.~\ref{triplepoint}.

\item
The Alexander numbering ${\cal L}$ of  regions divided by a given diagram 
is defined similarly.

\item
The source region $R$ and the
Alexander numbering $\LX(\tau)=\LX(R)$ 
are  defined for a triple point
$\tau$  using 
orientation 
normals.

\item
The sign $\epsilon(\tau)$ of a triple point $\tau$
is  defined \cite{CS:book} 
in such a way that it is positive if and only if 
the normals ot top, middle, bottom sheets, in this order, 
match the orientation of $3$-space. 

\item
For a coloring ${\cal C}$, the Boltzman weight 
at a triple point $\tau$ is defined by 
$B_T(\tau, {\cal C})=$ \linebreak 
$[\theta(x,y,z)^{\epsilon (\tau)} ]^{T^{-\LX(\tau)}}$, 
where $\theta$ is a $3$-cocycle, $\theta \in Z^3_{\rm TQ}(X;A)$. 
In the right of Fig.~\ref{triplepoint},
the triple point $\tau$ is positive, and ${\cal L}(\tau)=0$, 
so that $B_T(\tau, {\cal C})=\theta(p,q,r)$. 

\item 
The state-sum is defined by $\Phi_T(K)=  
\sum_{{\cal C}}  \prod_{\tau}  B_T( \tau, {\cal C}).
$

\end{itemize}

By checking the analogues of Reidemeister moves for knotted
surface diagrams, called Roseman moves, 
it was shown in \cite{CES} that 
$\Phi_T(K)$ is an invariant, 
called  the {\it (twisted quandle) cocycle invariant}
of knotted surfaces. 

Similarly, the state-sum invariant in the untwisted case was defined 
earlier 
in  \cite{CJKLS:era} and \cite{CJKLS}. In the untwisted case, 
there is no Alexander numbering, and the Boltzmann 
weight at  a triple point is  simply the quantity
$B(\tau, {\cal C})= \theta(x,y,z)^{\epsilon (\tau)} $ 
where  $x,y,z$ are the colors on the source regions of the bottom, middle,
and top sheets at the triple points.

In all of these cases, the value of the state-sum invariant depends
only on the cohomology class represented by the defining cocycle.
In particular, a coboundary will simply count the number of colorings
of a knot 
or knotted surface 
by the quandle $X$.

\subsection*{Applications}

Two important topological applications have been obtained
using the cocycle invariants.

\begin{itemize}
\item
The $2$-twist spun trefoil $K$ and its orientation-reversed counterpart $-K$
have shown to have distinct cocycle invariants
using a cocycle in $Z^3_Q(R_3;\Z_3)$, 
providing a proof that $K$ is non-invertible \cite{CJKLS}.

The higher genus surfaces obtained from $K$ by adding arbitrary
number of trivial $1$-handles are also non-invertible, since 
such handle additions do not alter the cocycle invariant.

We note, again, that this result in higher genus cases
is not immediately obtained from \cite{Hillman,Ruber}, although
higher genus generalizations of the Farber-Levine pairing
\cite{Kawa90a} can be used. 

\item
The projection of the $2$-twist spun trefoil was shown to 
have at least four triple points \cite{SatShima}. 

The same cocycle group  $Z^3_Q(R_3;\Z_3)$, 
but a different cocycle found in \cite{Mochi}
 was 
used. 

\end{itemize}

\section{Virtual Knots  and Quandle Homology} 

In this section, we describe $2$-dimensional quandle homology classes 
as cobordism classes of quandle colored virtual knot diagrams. 
See \cite{CKS:geo,FRS,Flower,Greene} for more general geometric 
descriptions of homology classes. 

\begin{figure}
\begin{center}
\mbox{
\epsfxsize=2in
\epsfbox{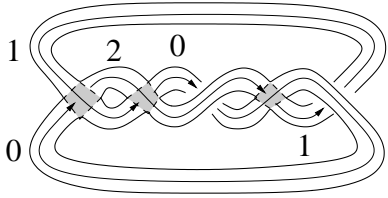} 
}
\end{center}
\caption{ A colored virtual knot }
\label{ssvsv0} 
\end{figure}

Consider an untwisted quandle homology class  of a 
quandle $X$, 
and represent the class by $\eta \in Z_2^{\rm Q}(X)$. 
Write $\eta$ as a sum of $2$-chains
$\eta = \sum_j \epsilon_j(a_j,b_j)$ where $\epsilon_j= \pm 1$. 
For each $j$  with $\epsilon_j=1$,
consider a positive crossing  diagram in which the over-arc is colored
$b_j$ and the under-arc away from which the normal to the over arc points
is colored $a_j$.  
Similarly, 
when $\epsilon_j=-1$ we consider a negative crossing of the 
same form. The boundary of the chain is $\pm (a_j - a_j*b_j)$ 
which is the difference in the colors on the under-arcs. 
Since $\eta$ is a cycle these boundary terms cancel over 
the sum of the crossings.

Thus to represent the $2$-cycle, we take a disjoint union 
of colored crossings, and join the end-point arcs together 
when they have the correct orientation and the same color. 
The arcs are joined together formally, and the joining need not
occur on a planar diagram,
obtaining a colored ``virtual knot diagram.''
Virtual knots have been popularized by L.H. 
Kauffman who has found, for example, that the diagram 
in Fig.~\ref{ssvsv0} has trivial Jones polynomial. 
A virtual knot can be regarded as a knot on a surface \cite{KK}.

Conversely, a colored virtual knot diagram represents a $2$-cycle. 
In Fig.~\ref{ssvsv0}, such a diagram colored with $R_3=\{0,1,2\}$
is depicted. The colored crossings in shaded squares,
from left to right,  represent $2$-chains 
$(0,1)$, $(1,2)$, and $-(1,0)$, respectively,
and therefore, this diagram shows that the $2$-chain 
$ (0,1) + (1,2) -(1,0)$ is a $2$-cocycle. 
The unshaded crossings between bands can be regarded as virtual crossings.
These bands connecting shaded squares correspond to identifying 
matching boundaries in the above construction.
Some remarks are in order.

\begin{itemize}
\item
There is a one-to-one correspondence  \cite{CKS:geo} between
(1) quandle colored virtual knot diagrams modulo the virtual Reidemeister moves
and  colored cobordisms, and   
(2) $2$-dimensional quandle homology classes.

\item
This geometric representations were used to estimate the rank of
rack homology groups in \cite{Flower,Greene} for some racks and quandles.

\item
This was also used to show that a certain long exact
sequence splits \cite{CKS:geo} 
at low dimensions, as mentioned in Section~\ref{quanss}.

\item 
The $2$-dimensional
 regions near crossings
in classical diagrams 
 can also be colored to represent $3$-cycles. 
Such colorings were used in \cite{RStre} 
to give an alternate proof that left- and right handed trefoils 
are not equivalent. 

\item 
The cocycle invariants can be interpreted as a 
formal sum of the Kronecker product between 
a fixed cocycle and such cycles constructed above represented by 
colored diagrams. Such an interpretation was used in \cite{CKS:sanfran}
to evaluate the cocycle invariants.

\item
Twisted cycles have a similar interpretation, but the 
consistency of Alexander numbering requires care.

\item 
The untwisted homology group $H_2^{\rm Q}(R_3)$ is trivial.
 Thus the cycle $(0,1)+(1,2)-(1,0)$ is a boundary. Meanwhile, 
the fundamental quandle \cite{LK:virt} of the virtual knot in
 Fig.~\ref{ssvsv0} can be computed to be $R_3$. Thus we have the interesting
 situtation in which a knot,
with any coloring by its fundamental quandle elements, 
 is null-homologous in the $2$-dimensional cycle group of  
 its fundamental quandle.

\end{itemize}

\section{Constructions of Cocycles from Extension Theory of Quandles}

The first constructions of quandle cocycles were a combination of hand  and computer calculations
\cite{CJKLS,CJKS1}. Here we summarize two important cases.
To describe cocycles, denote the characteristic function by 
$$\chi_{a}(x)= \left\{ 
\begin{array}{ll} 1 & {\mbox{\rm if }} \ a= x \\
0 & {\mbox{\rm if }} \ a \ne  x \end{array} \right. , $$
where $a$, $x$ are $n$-tuples of elements of a quandle $X$.

\begin{itemize} 
\item
For the Alexander quandle $S_4=  \Z_2[T,T^{-1}]/(T^2+T+1)$, 
the $\Z_2$-valued function
$$\phi =\sum_{a\ne b \;  a,b \ne T} \chi_{a,b}$$ represents a non-trivial 
cohomology class in $H^2_{\rm Q}(S_4; \Z_2) \cong \Z_2$.  

\item 
It was computed that $H^3_{\rm Q}(R_3;\Z_3)\cong \Z_3$ 
and 
a generator is 
given by 
$$\theta = -\chi_{(0,1,0)} + \chi_{(0,2,0)} -
\chi_{(0,2,1)}+\chi_{(1,0,1)}
    +\chi_{(1,0,2)}+\chi_{(2,0,2)}+\chi_{(2,1,2)}  
\in Z^3_{\rm Q} (R_3; \Z_3). $$

\end{itemize}

In \cite{CJKLS} 
it was mentioned that the trefoil ($3_1$) 
and the figure-eight knot $(4_1)$
have non-trivial cocycle invariants with the cocycle $\phi$. It was also shown that the $2$-twist spun trefoil is not invertible using the $3$-cocycle $\theta$.
This was proven using similar techniques in \cite{RStre}. 
Recently, Satoh and Shima  \cite{SatShima} have shown that any diagram
for the $2$-twist spun trefoil has at least 4 triple points using a 
$3$-cocycle 
in $Z^3_{\rm Q}(R_3; \Z_3)$ 
discovered by Mochizuki \cite{Mochi}.
Mochizuki  \cite{Mochi},  Litherland and Nelson \cite{LN} have developed more techniques for computing quandle homology and cohomology.

\bigskip

For quantum invariants,  solutions (R-matrices) to the Yang-Baxter equations 
were discovered by calculations first, and then Drinfeld
\cite{Drin} developed a theory of quantum groups whose 
representations gave rise to R-matrices. This construction is 
seen as an obstruction to  co-commutativity satisfying 
the next order (the Yang-Baxter) relation, or, deformation theory of
an algebraic structure giving rise to a solution to 
a 
higher order relation.
Considering analogies between group and quandle cohomology theories,
it is, then, natural to seek such methods of finding cocycles in 
deformation and extension theories of quandles. 
An extension theory of quandles was  
developed in \cite{CES}  for the 
twisted case as follows
(see also \cite{CENS,CKS:sanfran}), 
in analogy with the group cohomology theory
(one sees that the following is in parallel to Chapter IV of \cite{Brown}).

\begin{itemize}
\item
Let $X$ be a quandle and $A$ be an Alexander quandle,
so that $A$ admits an action by $\Z$ whose generator is denoted by $T$. Let $\phi \in Z^2_{\rm TQ} (X; A)$.
Let $AE(X,A, \phi)$ be the quandle defined on the set $A \times X$ by
the operation 
$(a_1, x_1) * (a_2, x_2) = (a_1 * a_2 + \phi(x_1, x_2), x_1 * x_2)$.

\item
The above defined operation $*$ on $A \times X$ indeed
defines a quandle $AE(X,A, \phi) = (A \times X, *)$,
which is called an {\em Alexander extension} 
of $X$ by $(A, \phi)$.

\item 
Let $X$ be a quandle and $A$ be an Alexander quandle.
Recall that $\eta \in Z^1_{\rm TQ}(X;A)$ implies that
$\eta: X \rightarrow A$ is a quandle homomorphism.
Let $0\rightarrow N  \stackrel{i}{\rightarrow} G
  \stackrel{p}{\rightarrow} A \rightarrow 0$ be an exact sequence
of $\Z [T, T^{-1}]$-module homomorphisms among
Alexander quandles. Let $s: A \rightarrow G$ be 
a set-theoretic section
(i.e., $ps=$id$_A$) with the ``normalization condition'' $s(0)=0$.
Then $s\eta: X \rightarrow G$ is a mapping, which is 
not necessarily a quandle homomorphism.
We measure the failure by $2$-cocycles. 
Since $p[T s\eta(x_1)+ (1-T) s\eta (x_2)]=p[ s\eta(x_1*x_2)]$
for any $x_1 , x_2 \in A$, there is $\phi(x_1, x_2) \in N$
such that 
$$
T s\eta(x_1)+  s\eta (x_2)= i \phi(x_1, x_2) + 
[ T  s\eta (x_2) +  s\eta(x_1*x_2)] . 
$$ 
This defines a function $\phi \in C^2_{\rm TQ}(X;N)$.
Then it was shown that 
$\phi \in Z^2_{\rm TQ}(X;N)$.

\item 
Let $s': A \rightarrow G$ be another section, 
and $\phi' \in  Z^2_{\rm TQ}(X;N)$
be a $2$-cocycle determined by 
$$ 
T s'\eta(x_1)+  s'\eta (x_2)= i \phi'(x_1, x_2) + 
[ T  s'\eta (x_2) +  s'\eta(x_1*x_2)] . 
$$ 

Then it was shown that 
$[\phi]=[\phi'] \in  H^2_{\rm TQ}(X;N)$.

\item 
It was shown that 
if $[\phi]=0 \in  H^2_{\rm TQ}(X;N)$,
then $\phi$ extends to a quandle homomorphism to $G$, 
i.e., there is a quandle homomorphism $\eta': X \rightarrow G$
such that $p \eta'=\eta$.

\end{itemize}

The above results were summarized as

\begin{theorem} {\bf \cite{CES}} \label{2cocyobstthm}
The obstruction  to extending $\eta: X \rightarrow A$ to 
a quandle homomorphism $X \rightarrow G$ lies
in  $H^2_{\rm TQ}(X;N)$.
\end{theorem}

Conversely, we have the following. 

\begin{lemma} {\bf \cite{CES}} \label{cocylemma}
Let $X$, $E$ be  quandles, and $A$ be an Alexander quandle.
Suppose there exists a bijection 
$f: E \rightarrow A \times X$ with the following property.
There exists a function $\phi: X \times X \rightarrow A$ such that
for any $e_i \in E$ ($i=1,2$), 
if $f(e_i)=(a_i, x_i)$, then 
$f(e_1 *  e_2) = (a_1 *a_2 + \phi(x_1, x_2) , x_1 * x_2 )$. 
Then $\phi \in Z^2_{\rm TQ}(X; A)$.
\end{lemma}

This lemma implies that under the same assumption  
we have $E=AE(X,A,\phi)$,  where  $\phi \in Z^2_{\rm TQ}(X; A)$,
and by identifying such quandles, we obtain cocycles as desired. 
We identify such examples, and include a proof, as it provides
explicit formulas of cocycles. 

Let $\Lambda_p=\Z_p[T, T^{-1}]$ for a positive integer $p$
(or $p=0$, in which case  $\Lambda_p$ is understood to be 
$\Lambda=\Z[T, T^{-1}]$). 
Note that  since $T$ is a unit in $\Lambda_p$,
$\Lambda_p/(h)$ for a Laurent polynomial $h \in \Lambda_p$ is
isomorphic to $\Lambda_p/(T^n h)$ for any integer $n$, so that 
we may assume that $h$ is a polynomial with a non-zero constant 
(without negative exponents of $T$).

\begin{lemma} {\bf \cite{CES}} \label{coeffcocylemma}
Let $h \in \Lambda_{p^{m}}$ be a polynomial with  leading and 
constant coefficients invertible, or $h=0$.
Let $\bar{h} \in  \Lambda_{p^{m-1}}$
and $\tilde{h} \in  \Lambda_{p}$ be such that 
$\bar{h} \equiv h \ \mbox{mod} \ (p^{m-1})$ 
and $\tilde{h} \equiv h \ \mbox{mod} \ (p)$, respectively
(in other words, $\bar{h}$ is $h$ with its
coefficients reduced modulo $p^{m-1}$, and $\tilde{h}$ is $h$
with its coefficients reduced modulo $p$).
Then the   quandle $E=\Lambda_{p^m}/({h})$
satisfies the conditions in Lemma~\ref{cocylemma}
with $X=\Lambda_{p^{m-1}}/(\bar{h})$ and $A=\Lambda_p / (\tilde{h})$. 

In particular, $\Lambda_{p^m}/({h}) $ is an Alexander extension
of $\Lambda_{p^{m-1}}/(\bar{h})$ by
$\Lambda_p / (\tilde{h})$:
$$\Lambda_{p^m}/(h)
= AE(\Lambda_{p^{m-1}}/(\bar{h}), \ \Lambda_p / (\tilde{h}) , \ \phi) , $$
for some 
$\phi \in Z^2_{\rm TQ}(\Lambda_{p^{m-1}}/(\bar{h}); \Lambda_p / (\tilde{h}) )$.
\end{lemma}
{\it Proof.\/} 
Let $A \in \Z _{p^m}$. Represent $A$ in $p^m$-ary notation
as $$A=\sum_{i=0}^{m-1} A_i p^i$$ where 
$A_i \in \{0, \ldots, p-1\}.$
Since $p$ is fixed throughout, we  represent $A$ by the sequence
$$[A_{m-1}, A_{m-2}, A_{m-3}, \ldots, A_0 ]. $$ 
Define $\overline{A} = [A_{m-2}, \ldots, A_0].$
Observe that $A \equiv \overline{A} \pmod{ p^{m-1}}$, and
$A \equiv A_0 \pmod{p}$. 

Let  $\hat{\pi}: \Z_{p^m} \rightarrow \Z_{p^{m-1}}$ be the map 
defined by $\hat{\pi}(A) = \overline{A}$. We obtain a short exact
sequence:

$$0 \rightarrow \Z_p \stackrel{\hat{\imath}}{\rightarrow} \Z_{p^m}
\stackrel{\hat{\pi}}{\rightarrow} \Z_{p^{m-1}} \rightarrow 0$$
where $\hat{\imath}(A)=[A,0,\dots,0]$. There is a set-theoretic section
$ \Z_{p^m} \stackrel{\hat{s}}{\leftarrow} \Z_{p^{m-1}}$ defined by
$\hat{s}[A_{m-2},\dots,A_0]=[0,A_{m-2},\dots,A_0].$
The map $\hat{s}$ satisfies $\hat{\pi}\hat{s}={\rm id}$ and 
$\hat{s}(0)=0$.

For a polynomial $L(T) \in \Lambda_{p^m} = \Z_{p^m}[T, T^{-1}]$,
write $$L(T) = \sum_{j=-n}^k [A_{j,m-1}, A_{j, m-2}, \ldots, A_{j,0}]
T^j.$$
Define $$\overline{L}(T) = \sum_{j=-n}^k [ A_{j, m-2}, \ldots, A_{j,0}] T^j \in
\Lambda_{p^{m-1}},$$ 
and 
$$\tilde{L}(T) = \sum_{j=-n}^k A_{j,m-1} T^j \in \Lambda_{p}.$$
There is a one-to-one correspondence $f: \Lambda_{p^m} \rightarrow
\Lambda_{p} \times \Lambda_{p^{m-1}}$ given by
$f(L)= (\tilde{L}, \overline{L})$. 
We have a short exact sequence of rings:

$$0 \rightarrow \Z_p[T,T^{-1}] \stackrel{i}{\rightarrow}
\Z_{p^m}[T,T^{-1}] \stackrel{\pi}{\rightarrow} \Z_{p^{m-1}}[T,T^{-1}]
\rightarrow 0$$
with a set theoretic section
$ \Z_{p^m}[T,T^{-1}] \stackrel{s}{\leftarrow} \Z_{p^{m-1}}[T,T^{-1}]$
where $i$, $\pi$ and $s$ are the natural maps induced by 
$\hat{i}$, $\hat{\pi}$ and $\hat{s}$, respectively.
Note that 
for $L \in \Lambda_{p^m} = \Z_{p^m}[T,T^{-1}]$ we have
$\overline{L}= \pi(L)$, and 
the section $s:\Lambda_{p^{m-1}} \rightarrow \Lambda_{p^m}$ is defined by
the formula
$$s \left(\sum_{j=-n}^k [ A_{j, m-2}, \ldots, A_{j,0}] T^j\right) = 
\sum_{j=-n}^k [0, A_{j, m-2}, \ldots, A_{j,0}] T^j.$$

For $L, M \in \Lambda_{p^m}$, let  
$$s(\overline{L})*s(\overline{M})= \sum_{j} [F_{j, m-1}, \ldots, F_{j, 0}] T^j \in \Lambda_{p^{m-1}}.$$
If $L= \sum_j A_j T^j$, and $M= \sum_j B_j T^j$, then 
$$L*M =  B_{-n} T^{-n} + \sum_{j=-n+1}^{k+1} (A_{j-1} - B_{j-1} +B_j)
T^j = \sum_{j=-n}^k C_j T^j.$$
Furthermore, 
\begin{eqnarray*}
\overline{L} *\overline{M}  &=& [B_{-n,m-2}, \ldots, B_{-n,0}] T^{-n} \\
           &+& \sum_{j=-n+1}^{k+1} \left( [A_{j-1,m-2}, \ldots ,
A_{j-1,0}] -[B_{j-1,m-2}, \ldots , B_{j-1,0}] +[B_{j,m-2}, \ldots ,
B_{j,0}] \right) T^j
\end{eqnarray*} 
and write the right-hand side by $ \sum_{j=-n}^k D_j T^j$.
Note that $D_j$'s are well-defined integers, 
not only  elements of $\Z_{p^{m-2}}$.
If $D_j$ is positive, then $F_{j, m-1}=0$, and if 
$D_j$ is negative, then $F_{j, m-1}=p-1$.
Hence
$$f(L*M)=(\tilde{L}*\tilde{M} +
\phi(\overline{L} , \overline{M}), \overline{L}*\overline{M}), $$
where 
$$ \phi( \overline{L}, \overline{M})=\sum_{j=-n}^k F_{j, m-1}. $$
This concludes the case $h=0$.

Now let $h(T) \in \Z_{p^m}[T]$  be a  polynomial with  leading and constant
coefficients being 
invertible in $\Z_p$. 
  Let $(h)$ denote the ideal generated by $h$.
Since  $i(\tilde{h}) \subset (h)$,  
we obtain a short exact
sequence of quotients:

$$0 \rightarrow \Z_p[T,T^{-1}]/(\tilde{h}) \stackrel{\overline{\imath}}{\rightarrow}
\Z_{p^m}[T,T^{-1}]/(h) \stackrel{\overline{\pi}}{\rightarrow}
\Z_{p^{m-1}}[T,T^{-1}]/(\overline{h}) \rightarrow 0$$
with a set-theoretic section 
$ \Z_{p^m}[T,T^{-1}]/(h) \stackrel{\overline{s}}{\leftarrow}
\Z_{p^{m-1}}[T,T^{-1}]/(\overline{h}).$
Thus we obtain a twisted cocycle 
$$\phi: \Z_{p^{m-1}}[T,T^{-1}]/(\overline{h}) \times
\Z_{p^{m-1}}[T,T^{-1}]/(\overline{h}) \rightarrow \Z_p[T,T^{-1}]/(\tilde{h})
\quad \qed $$

Since $R_n=\Lambda_n / (T+1)$, we have the following.

\begin{corollary} \label{dihedcocylemma}
The dihedral quandle $E=R_{p^m}$, where $p,m$ are positive integers with $m>1$,
satisfies the conditions in Lemma~\ref{cocylemma}
with $X=R_{p^{m-1}}$ and $A=R_p$. 

In particular, $R_{p^m} $ is an Alexander extension
of $R_{p^{m-1}}$ by $R_p$: $R_{p^m}= AE(R_{p^{m-1}}, \  R_p, \ \phi)$, 
for some $\phi \in Z^2_{\rm TQ}(R_{p^{m-1}}; R_p)$.
\end{corollary}

\begin{example} \label{dihedex} {\rm
Let $X=R_3$ and $A=R_3$, then the proof of Lemma~\ref{coeffcocylemma}
gives an explicit $2$-cocycle $\phi$ as follows. For 
$\phi(r_1, r_2)=\phi(1,2)$, for example, one computes
$$ r_1 * r_2 = [0,1]* [0,2]= 2 [0,2]- [0,1] =3 =3 \cdot 1 + 0=[1,0],$$ 
Hence $\phi(0,2)=1$. 
In terms of  the characteristic function, 
the cocycle $\phi$ contains the term $\chi_{0,2}$.
By computing the quotients for all pairs, one obtains
$$ \phi = \chi_{0,2} + \chi_{1,2} +  2  \chi_{1,0} + 2 \chi_{2,0}. $$

The same argument was applied to $R_{\infty} $
to show  that the quandle $R_{\infty} $ is an Alexander extension of $R_n$ by 
$R_{\infty} $, for any positive integer $n$.

} \end{example}

Similar techniques give us untwisted cocycles 
\cite{CENS}, with  explicit formulas for these $2$-cocycles
as follows. 
In this case, the extension is called an {\it abelian extension},
denoted by $E=E(X, A, \phi) $ for $\phi \in Z^2_{\rm Q}(X;A)$, 
and the quandle operation on $E=A \times X$ is defined  by 
$(a_1, x_1)*(a_2, x_2)=(a_1 + \phi(x_1, x_2), x_1 * x_2)$. 

\begin{itemize}
\item
For any  positive integers $q$ and $m$, 
$E=\Z _{q^{m+1}} [T, T^{-1}] /  (T -1 +q) $ is an abelian extension 
$E=E( \Z _{q^{m}} [T, T^{-1}] /  (T -1 +q) ,  \Z_q,  \phi)$
of $X= \Z _{q^{m}} [T, T^{-1}] /  (T -1 +q)$ for some cocycle
$\phi \in  Z^2_{\rm Q}( X;  \Z_q)$.

\item 
For any positive integer $q$ and $m$, the quandle 
$E=\Z_q [T, T^{-1} ] / (1-T)^{m+1} $ is an abelian extension
of $X=\Z_q [T, T^{-1} ] / (1-T)^{m} $ over $\Z_q$: 
$E=E(X, \Z_q, \phi)$, for some $\phi \in Z^2_{\rm Q}(X; \Z_q)$.

\end{itemize}

Furthermore, for untwisted $2$-cocycles, 
an interpretation of the cocycle knot invariant 
was given  \cite{CENS}
as an obstruction
to extending a given coloring of a knot diagram by a quandle $X$ 
to a coloring by an abelian extension $E$. 
Similar interpretations for twisted case or knotted surface case are 
unknown.

Ohtsuki \cite{Ohtsuki}
defined a new  cohomology theory for quandles and 
an 
extension theory, 
together with a list of problems in the subject. 

\subsection*{Acknowledgements}

We gratefully acknowledge the contributions of our collaborators in these projects:
Mohammed Elhamdadi, Daniel Jelsovsky, Seiichi Kamada, Louis Kauffman, Laurel Langford,  and Marina Nikiforou. We thank the organizers of the 2001 GTC for their hard work and for allowing us the opportunity to present these results. As we were finishing this paper we received a copy of the survey \cite{Kam:sur}. We are pleased to refer the reader to that paper as well.

\end{document}